\documentclass[12pt]{amsart}
\usepackage{amsmath,amssymb,amsthm}
\usepackage{color,graphics}
\usepackage{srcltx}

\makeatletter

\newtheorem{theorem}{Theorem}[section]
\newtheorem{lemma}[theorem]{Lemma}
\newtheorem{corollary}[theorem]{Corollary}
\newtheorem{remark}[theorem]{Remark}

\theoremstyle{definition}

\newtheorem{example}[theorem]{Example}

\newtheorem{problem}[theorem]{Problem}

\makeatother

\def\per{\mathop{\mathrm{per}}\nolimits}
\def\Id{\mathop{\mathrm{Id}}\nolimits}

\def\tr{{\rm tr}\,}
\def\FF{{\mathbb F}}

\def\CC{{\mathbb C}}

\def\charac{\mathop{\mathrm{ch}}}

\def\rk{{\rm rk \,}}

\def\st{\widehat}

\newcounter{myenumi}

\newsavebox{\cmm}
\savebox{\cmm}{\indent}
\newenvironment{myenumerate}[1]{
\begin{list}{
{\bf #1~\themyenumi}. } {\labelwidth=0pt
\labelsep=0pt\leftmargin=0pt\usecounter{myenumi}} }{\end{list}}

\begin{document}
\thispagestyle{empty}
\begin{center}
{\bf \Large On the Polya permanent problem over finite fields}
\end{center}\vspace{2cm}
Gregor Dolinar\\
 Faculty of Electrical Engineering, University of Ljubljana,
Tr\v{z}a\v{s}ka 25, SI-1000 Ljubljana, Slovenia.\\

\noindent{Alexander E.\ Guterman}\\ Faculty of Algebra, Department of Mathematics and
Mechanics,
Moscow State University, GSP-1, 119991 Moscow, Russia.\\

\noindent{Bojan Kuzma}\\ {${}^1$University of Primorska, Glagolja\v{s}ka 8, SI-6000
Koper, Slovenia,
\and ${}^2$IMFM, Jadranska 19, SI-1000 Ljubljana, Slovenia.}\\

\noindent{Marko Orel}\\ {IMFM, Jadranska 19, SI-1000 Ljubljana, Slovenia.}
\newpage
\thispagestyle{empty}
\noindent{\bf Proposed running head:} Polya permanent problem\\

\noindent{\bf Corresponding author:}\\
 {Bojan Kuzma} IMFM, Jadranska 19, SI-1000 Ljubljana,
Slovenia.\\
 {\bf e-mail: \verb"bojan.kuzma@pef.upr.si"}
\newpage
\addtocounter{page}{-2}
\title[Polya permanent problem]{On the Polya permanent problem over finite fields}
\author[Dolinar]{Gregor Dolinar}
\address[Gregor Dolinar]{Faculty of Electrical Engineering, University of Ljubljana,
Tr\v{z}a\v{s}ka 25, SI-1000 Ljubljana, Slovenia.}
 \email[Gregor Dolinar]{gregor.dolinar@fe.uni-lj.si}
 \author[Guterman]{Alexander E.\ Guterman}
\address[Alexander E.\ Guterman]{Faculty of Algebra, Department of Mathematics and
Mechanics, Moscow State University, GSP-1, 119991 Moscow, Russia.}
 \email[Alexander E.\ Guterman]{guterman@list.ru}
\author[Kuzma]{Bojan Kuzma}
\address[Bojan Kuzma]{${}^1$University of Primorska,
Glagolja\v{s}ka 8, SI-6000 Koper, Slovenia, \and ${}^2$IMFM, Jadranska 19, SI-1000
Ljubljana, Slovenia.}
 \email[Bojan Kuzma]{bojan.kuzma@pef.upr.si}
 \author[Orel]{Marko Orel}
\address[Marko Orel]{IMFM, Jadranska 19, SI-1000 Ljubljana, Slovenia.}
 \email[Marko Orel]{marko.orel@fmf.uni-lj.si}


\thanks{The research was supported by a joint Slovene--Russian grant
BI-RU/08-09-009. The research of the second author is also supported by the RFBR grant
08-01-00693a}

\begin{abstract}
Let $\FF$ be a finite field of characteristics different from two. We show that no
bijective map transforms permanent into determinant when the cardinality of $\FF$ is
sufficiently large. We also give an example of non-bijective map when $\FF$ is arbitrary
and  an example of a bijective map  when $\FF$ is infinite which do transform permanent
into determinant. The developed technique  allows us to estimate the probability of the
permanent and the determinant of matrices over finite fields to have a given value. Our
results are also true over finite rings without zero divisors.
\end{abstract}

 \maketitle

\section{Introduction}
Let $A=(a_{ij})\in M_n(\FF)$ be an $n\times n$ matrix over a field $\FF$. The permanent
function
$$ \per A= \sum_{\sigma\in S_n} a_{1\sigma(1)}\cdots a_{n\sigma(n)} $$
is defined in a very similar way to the definition of the determinant function
$$ \det A= \sum_{\sigma\in S_n} sgn({\sigma}) a_{1\sigma(1)}\cdots a_{n\sigma(n)} .$$
In both cases the sum is considered over all permutations $\sigma\in S_n$, where $S_n$
denotes the set of all permutations of the set $\{1,2,\ldots, n\}$. The value
$sgn(\sigma)\in \{-1,1\}$ is the signum of the permutation $\sigma$, i.e.,
$sgn(\sigma)=1$ if $\sigma$ is an even permutation, and $sgn(\sigma)=-1$ if $\sigma$ is
an odd permutation.

The determinant  is certainly one of the most well-studied functions in mathematics.
Geometrically, it is the volume together with orientation of the parallelepiped defined
by rows (or columns) and algebraically, it is the product of all eigenvalues, counted
with their multiplicities. The permanent function is also well-studied, especially in
combinatorics, see~\cite{mink}. For example, if $A$ is a (0,1)-matrix, then the value
$\per A$ is equal to the number of perfect matchings in a bipartite graph with adjacency
matrix $A$. However, no nice geometric or algebraic interpretation is known for
permanent. Moreover the permanent does not enjoy the same properties as the determinant,
in particular it is neither multiplicative nor invariant under linear combinations of
rows or columns.

Computing permanent of a matrix seems to have different computational complexity than
computing the determinant. The determinant can be calculated by a polynomial time
algorithm. For example, Gauss elimination method requires $O(n^3)$ operations. At the
same time no efficient algorithm for computing the permanent function is known, and, in
fact, none is believed to exist.  When using its definition, the computation of the
permanent requires $(n-1)n!$ multiplications and one of the best known algorithms to
compute permanent, so-called Ryser's formula~\cite{ryser}, has an exponential complexity
and requires $(n-1)\cdot (2^n-1)$ multiplications. Moreover, Valiant~\cite{valiant} has
shown that even computing the permanent of a (0,1)-matrix is a $\sharp P$-complete
problem, i.e., this problem is an arithmetic analogue of Cook's hypothesis $P\ne NP$,
see~\cite{cook,karp,garey_johnson} for details.

Starting from  1913 researchers are trying to find a way to calculate permanents using
determinants. More precisely,  the following problems  which dates back to the work of
P\'olya~\cite{polya}  are under intensive investigations for almost a century.
\begin{problem} \label{pr_pm}
Does there exist  a uniform way of affixing $\pm$ signs to the entries of a matrix
$A=(a_{ij})\in M_n(\FF) $ such that $\per (a_{ij})=\det (\pm a_{ij})? $
\end{problem}
\begin{problem} \label{pr_even_cycles}
Given a (0,1)-matrix $A\in M_n(\FF)$, does there exists a transformed matrix $B$,
obtained by changing some of the $+1$ entries of $A$ into $-1$, so that $\per A=\det B$?
\end{problem}
\begin{problem} \label{pr_general}
Under what conditions does there exist a transformation $\Phi: M_n(\FF) \to  M_m(\FF)$
satisfying
\begin{equation}
\label{eq:per-det}
\per A=\det \Phi(A)?
\end{equation}
\end{problem}
In this case the image $\Phi(A)$ is usually called a {\em P\'olya matrix\/} for~$A$.

For example if $n=2$, one can consider the P\'olya matrix
\begin{equation}\label{eq:n=2} B=B(A)=\left( \begin{array}{cc}a_{11} & -a_{12} \\ a_{21} & a_{22} \end{array}\right).
\end{equation}

Problem~\ref{pr_pm} was solved negatively by Szeg\"o in \cite{szego}, namely he proved
that for $n\ge 3$ there is no generalization of the formula~(\ref{eq:n=2}).

Problem~\ref{pr_even_cycles} has been intensively studied since it belongs to the famous
class of equivalent problems, containing the following ones:

\noindent When does a real square matrix have the property that every real matrix with
the same sign pattern is non-singular? When does a bipartite graph have a ``Pfaffian
orientation''? Given a digraph, does it have no direct circuit of even length?
See~\cite{brualdi_shader,robertson_seymour_thomas,vazirani_yannakakis} for the detailed
and self-contained information.

Problem~\ref{pr_general} is  a natural generalization of Problem~\ref{pr_pm}. Namely,
affixing ($\pm1$) signs to the entries of a matrix is an example of a certain linear
transformation with easy structure. One may ask, if there exists some more sophisticated
linear transformation $\Phi:M_n(\FF) \to M_n(\FF)$ satisfying~(\ref{eq:per-det})? In 1961
Marcus and Minc~\cite{marcus_mink}, see also Botta~\cite{botta}, proved that if $n\ge 3$
there is no linear transformations $\Phi:M_n(\FF) \to M_n(\FF)$ satisfying
equality~(\ref{eq:per-det}).

Von zur Gathen~\cite{gathen} investigated the linear transformations $\Phi:M_n(\FF) \to
M_m(\FF)$ satisfying the equality
\begin{equation}
\label{eq:det-per}
\det A=\per \Phi(A)
\end{equation}
 and proved that if there exists such $\Phi$, then $m>\sqrt 2 n - 6 \sqrt n$. These results were later improved, see for example Cai~\cite{cai} and references therein.

After that several attempts to further reduce  the linearity assumption were made, see
for example Coelho, Duffner~\cite{coelho_duffner}, Kuzma~\cite{kuzma}, and references
therein. In these works no  bijectivity or linearity is assumed, but the authors consider
transformations $\Phi:M_n(\CC) \to M_n(\CC)$ satisfying the equality
 $$d_{\chi}(\Phi(A)+\lambda\Phi(B))=d_{\chi'}(A+\lambda B), \qquad (\lambda\in\CC)$$
 where $ d_{\chi},d_{\chi'}$ are arbitrary immanants.
In particular, this also covers the possibility $d_{\chi}=\det$ and $d_{\chi'}=\per$ or
vise versa.

In the present paper we show that without any regularity assumptions imposed on $\Phi$
there do exist  transformations  (possibly {\em bijective\/} if the underlying field is
infinite) that even exchange the permanent and the determinant, i.e. transformations that
satisfy both equalities~(\ref{eq:per-det}) and~(\ref{eq:det-per}) simultaneously, see
Examples~\ref{Ex:main} and~\ref{Ex:2} from the present paper.

The main aim of the present paper, however, is to show that if $\FF$ is a finite field of
sufficiently large cardinality, depending on $n$, then there are no bijective
transformations  $\Phi:M_n(\FF) \to M_n(\FF)$ satisfying~(\ref{eq:per-det}), i.e., we
obtain a negative solution of Problem~\ref{pr_general} for bijective maps, defined on
matrices over finite fields.

These results are heavily based on the detailed analysis of the cardinality of the set of
matrices over finite fields with zero permanent.

As an application of our results we also estimate the probability of the permanent and
the determinant of matrices over a finite field to have a given value. This problem dates
back to the works by Erd\"os and R\'enyi~\cite{erdos_renyi,erdos_renyi1}, where they
estimated the probability for a (0,1)-matrix with a given number of ones to have a zero
permanent. For the detailed and self-contained account of the results one may study
monograph~\cite{borovskikh_korolyuk,sachkov_tarakanov}.

Our paper is organized as follows:  Section~\ref{sec-statement} contains basic
definitions and notations used in the paper and the statements of the main results. In
Section~\ref{Sec:2} we calculate the number of matrices with the zero permanent in
$M_3(\FF)$. In Section~\ref{Sec:3} we compute  the cardinality of a set of pairs of
vectors, which are orthogonal to each other with respect to a matrix of fixed rank $r$
and split the set of all matrices with zero permanent into several subsets by means of
the Laplace decomposition. In Section~\ref{Sec:4} we introduce inductively the lower and
upper bounds for the number of matrices with zero permanent and prove these bounds.
Section~\ref{sec-proof} is devoted to the proof of the main result, i.e., that the
introduced upper bound for matrices with zero permanent is strictly less than the number
of matrices with zero determinant. In Section~\ref{Sec:7} the probability of the
permanent and the determinant to have given values is estimated. In Section~\ref{Sec:6}
we provide some examples, in particular the examples of non-bijective transformations on
matrices over any fields and bijective transformations of matrices over  infinite field
that satisfy both equalities~(\ref{eq:per-det}) and~(\ref{eq:det-per}) simultaneously. We
also show that our results are valid over finite rings without zero divisors.

\section{Preliminaries and statement of the main result}\label{sec-statement}

In our paper $\FF$ is a finite field of characteristics different from 2 and of the
cardinality $|\FF|=q$, except in Section~\ref{Sec:6}, where $\FF$ is arbitrary.

We denote the identity matrix from $M_n(\FF)$ by $I_n$ and zero matrix by $O_n$. If the
size $n$ is clear from the context, we omit the corresponding index. By $A^{\tr}$ we
denote the transposed matrix to $A$.

In this paper we use the term {\em per-minor of order $k$\/} (or {\em $k$-by-$k$
per-minor\/}) of $A\in M_n(\FF)$ to denote the permanent of a $k\times k$-submatrix of
$A$. {\em Principal\/} per-minor of $A\in M_n(\FF)$ is a per-minor of order $n-1$. Let
$A_{ij}$ denote the matrix obtained from $A$ by deleting the $i$-th row and $j$-th
column; $A_{(i\ldots j)(k\ldots l)}$ denote the matrix obtained from $A$ by deleting the
rows from $i$ to $k$ and the columns from $j$ to $l$. Let
\begin{equation*}\label{eq:hatA} \st{A}=\begin{pmatrix}
\per A_{11}& \dots & \per A_{1n}\\
\vdots &\ddots & \vdots\\
\per A_{n1} &\dots & \per A_{nn}
\end{pmatrix}
\end{equation*}
be a {\em permanental compound\/} of $A$. In this paper we investigate the sets
$$P_n(\FF) = \{A\in M_n(\FF):\;\; \per A=0\}$$ and $$D_n(\FF) = \{A\in M_n(\FF):\;\; \det
A=0\}$$ of all matrices with zero permanent and zero determinant, respectively.

It is straightforward to see that as is the case with a classical determinant, the
permanent  also obeys the Laplace decomposition, see for example~\cite[Chapter 2, Theorem
1, 2]{mink},
$$ \per A=\per (a_{ij})=a_{i1}\per A_{i1}+a_{i2} \per A_{i2}+\ldots+ a_{in} \per A_{in}.$$

\bigskip

 Our main result can be formulated as follows.
\begin{theorem}\label{theorem}
Suppose $n\ge 3$. Then there exists $q_0$, depending on $n$, such that for any finite
field $\FF$ with at least $q_0$ elements and  $\charac\FF\ne2$ no bijective map
$\Phi:M_n(\FF)\to M_n(\FF)$ satisfies
\begin{equation}\label{eq:main}
\per A=\det\Phi(A).
\end{equation}
When $n=3$ the conclusion holds for any finite field with $\charac\FF\ne2$.
\end{theorem}

 Note that any finite ring without zero
divisors is a field (in Section~\ref{Sec:6} we provide a short proof of this fact for the
sake of completeness). Therefore the above result is valid also for matrices over finite
rings without zero divisors.
\begin{corollary} \label{cor_rings}
Let $n\ge 3$ and let $R$ be a finite ring without zero divisors of sufficiently large
cardinality, $\charac R\ne2$. Then no bijective map $\Phi : M_n(R)\to M_n(R)$ satisfies
$\per(A)=\det \Phi(A)$.
\end{corollary}

\begin{remark}
By considering $\Psi=\Phi^{-1}$ the above Corollary shows that $\per\Psi(A)=\det(A)$ is impossible
for bijective $\Psi$ acting on matrices over a finite ring without zero divisors of sufficiently large
cardinality and characteristic different from~2.
\end{remark}

The proof of Theorem~\ref{theorem} will be given in Section~\ref{sec-proof}. Here we
outline the main idea. Any bijective~$\Phi$ satisfying~(\ref{eq:main}) would induce a
bijection from the set $P_n(\FF)$ of $n$--by--$n$ matrices with zero permanent onto the
set $D_n(\FF)$ of $n$--by--$n$ matrices with zero determinant. Consequently, to prove the
theorem it suffices to show that the number $|P_n(\FF)|$ does not equal $|D_n(\FF)|$ for
all sufficiently large finite fields of characteristic different from two. We remark that
the latter number is well-known. Actually, there exists precisely
\begin{equation}\label{eq:det=0}
|D_n(\FF)|=q^{n^2}-\prod_{k=1}^n(q^n-q^{k-1})= q^{n^2}-q^{n(n-1)\over
2}(q^n-1)\dots(q-1)
\end{equation} $n$--by--$n$ matrices with determinant zero~\cite[Prop.~2,
p.~41]{alperin_rowen}, where $q=|\FF|$.

For $n=3$ we will exactly calculate the number $|P_3(\FF)|$ of matrices with permanent
zero, however for $n\ge 4$ we will not give an exact formula for $|P_n(\FF)|$, but we
will give its upper bound $\mathfrak{U}_n(\FF)$ and show that
$$|P_n(\FF)|\le \mathfrak{U}_{n}(\FF)\lneq |D_n(\FF)|$$
if $\FF$ is a finite field with sufficiently many elements and $\charac\FF\ne2$.

\section{Zero permanents in $M_3(\FF)$ \label{Sec:2}}
\begin{lemma}\label{lem:permanent_for_3-by-3}
Let $\FF$ be a finite field  with $\charac\FF\ne2$. Then
$$|P_3(\FF)|=|D_3(\FF)|-q^2(q-1)^5.$$
\end{lemma}
\begin{proof}
We decompose $D_3(\FF)$ and $P_3(\FF)$ into pairwise disjoint union of three sets
$$D_3(\FF)=D_3'(\FF)\cup D_3''(\FF)\cup D_3'''(\FF)$$
and
$$P_3(\FF)=P_3'(\FF)\cup P_3''(\FF)\cup P_3'''(\FF),$$ where
\begin{align}
\nonumber D_3'(\FF)&=\{A\in D_3(\FF):\;\;a_{33}\ne0,\,\det A_{11}=0\}\\ 
D_3''(\FF)&=\{A\in D_3(\FF):\;\;a_{33}\ne0,\,\det A_{11}\ne0\}\label{eq:2}\\
D_3'''(\FF)&=\{A\in D_3(\FF):\;\;a_{33}=0\}\label{eq:3}\\[2mm]
\nonumber P_3'(\FF)&=\{A\in P_3(\FF):\;\;a_{33}\ne0,\,\per A_{11}=0\}\\
P_3''(\FF)&=\{A\in P_3(\FF):\;\;a_{33}\ne0,\,\per A_{11}\ne0\}\tag{\ref{eq:2}'}\\
P_3'''(\FF)&=\{A\in P_3(\FF):\;\;a_{33}=0\}\tag{\ref{eq:3}'}
\end{align}
We claim that the cardinality of (\ref{eq:2}) and (\ref{eq:2}') are the same; and also
the cardinality of (\ref{eq:3}) and (\ref{eq:3}') are the same, while $|P_3'(\FF)|=
|D_3'(\FF)|-q^2(q-1)^5$.

Start with (\ref{eq:2}).  There are $q^2(q-1)^2$ many ways of prescribing the values to
`variables' $a_{22},a_{23},a_{32},a_{33}$ to achieve $\det A_{11}\ne0\ne a_{33}$. We can
further arbitrarily prescribe the values of $a_{12},a_{13},a_{21},a_{31}$, while $a_{11}$
is then completely determined by $\det A=0$, i.e., by
$$a_{11}=
   \frac{a_{13} \left(a_{22}
   a_{31}-a_{21}
   a_{32}\right)+a_{12}
   \left(a_{21}
   a_{33}-a_{23}
   a_{31}\right)}{a_{22}
   a_{33}-a_{23}
   a_{32}}.$$
In total, $|D_3''(\FF)|=q^2(q-1)^2\cdot q^4$. A similar computation also gives
$|P_3''(\FF)|=q^2(q-1)^2\cdot q^4$, as claimed.

We next show that the cardinalities of (\ref{eq:3}) and (\ref{eq:3}') are the same. To do
this, just notice that
$$\Psi_{33}:\left(
\begin{array}{lll}
 x_{11} & x_{12} & x_{13}
   \\
 x_{21} & x_{22} & x_{23}
   \\
 x_{31} & x_{32} & 0
\end{array}
\right)\mapsto \left(
\begin{array}{lll}
 -x_{11} & x_{12} & x_{13}
   \\
 x_{21} & -x_{22} & x_{23}
   \\
 x_{31} & x_{32} & 0
\end{array}\right)$$ is a linear bijection with the property $\per
A=\det\Psi_{33}(A)$ for every $A \in M_3(\FF)$ with $a_{33} = 0$. Whence it also maps the
set $P_3'''(\FF)$ bijectively onto the set $D_3'''(\FF)$, as claimed.

Finally, we compute the cardinalities $|P_3'(\FF)|$ and $|D_3'(\FF)|$. Consider first the
set
\begin{align*}
D_3'(\FF)&=\{A\in D_3(\FF):\;\;a_{33}\ne0,\,\det A_{11}=0\}\\
 &= \Bigl\{\bigl(a_{ij}\bigr)\in M_3 (\FF):\;\;\tfrac{\left(a_{13}
   a_{32}-a_{12}
   a_{33}\right)
   \left(a_{21}
   a_{33}-a_{23}
   a_{31}\right)}{a_{33}}=0,\\
   &\hspace{6.0cm}\mbox{}a_{22}-
   \tfrac{a_{23}
   a_{32}}{a_{33}}=0,\, a_{33}\ne0\Bigr\}.
\end{align*}
The number of matrices inside $D_3'(\FF)$ can be computed as follows. We can choose a
total of $q-1$ distinct nonzero values for $a_{33}$, and a total of $q$  distinct values
for each `variable' $a_{23}$ and $a_{32}$. Once these are chosen, $a_{22}$ is uniquely
determined from them, by the second equation.  All together, we can prescribe $(q-1)q^2$
different values for `variables' $a_{33},a_{32},a_{23},a_{22}$.

Once we choose the values of these four `variables',  we have additional equation
${\left(a_{13}
   a_{32}-a_{12}
   a_{33}\right)
   \left(a_{21}
   a_{33}-a_{23}
   a_{31}\right)}=0$ with  four new `variables' $a_{13},a_{12},a_{21},
 a_{31}$. There are $q^4$ ways of prescribing their values, but only
 $q(q-1)\cdot q(q-1)$ ways of prescribing their values so that both
factors are nonzero --- in fact, we can prescribe, say, $a_{13}$ arbitrarily, and then
$a_{12}$ is determined by $a_{12}\ne a_{13}\frac{a_{32}}{a_{33}}$, likewise for the
second factor.  From this we deduce that there are precisely $q^4-q^2(q-1)^2=q^2(2q-1)$
ways of prescribing the values for $a_{13},a_{12},a_{21},
 a_{31}$ so that the product of the two factors is zero. Finally, due
 to $\det A_{11}=0$ we may arbitrarily prescribe $a_{11}$ without affecting
 $\det A=0$. In total,
$$|D_3'(\FF)|=(q-1)q^2\cdot q^2(2q-1)\cdot q=q^5(q-1)(2q-1).$$

In the set
\begin{align}
P_3'(\FF)&=\{A\in P_3(\FF):\;\;a_{33}\ne0,\,\per A_{11}=0\}\notag\\[2mm]
 &\label{eq:P3'}
 \begin{aligned}
 =\Bigl\{&\bigl(a_{ij}\bigr)\in M_3(\FF):\;\;\tfrac{a_{13} a_{32}
   \left(a_{21}
   a_{33}-a_{23}
   a_{31}\right)}{a_{33}}+\\
      &\mbox{}+a_{12} \left(a_{23}
   a_{31}+a_{21}
   a_{33}\right)=0,\,a_{22}+\tfrac{a_{23}
   a_{32}}{a_{33}}=0,\, a_{33}\ne0\Bigr\}
\end{aligned}
\end{align}
the first equation does not split, so we need a different approach to
compute~$|P_3'(\FF)|$. First we count those matrices inside~$P_3'(\FF)$ which satisfy
$\left(a_{23}
   a_{31}+a_{21}
   a_{33}\right)\neq 0$. As in the determinant case,  $a_{11}$ is arbitrary  while the values for `variables' $a_{33},a_{32},a_{23},a_{22}$ can be prescribed in
$(q-1)q^2$ different ways. Once these values are chosen, we have $q(q-1)$ possibilities
for `variables' $a_{31},a_{21}$ to achieve that $\left(a_{23}a_{31}+a_{21}
   a_{33}\right)\ne 0$. We may further prescribe $a_{13}$ arbitrarily, and
   then  $a_{12}$ is completely determined by
   $$a_{12}=   \frac{a_{13} a_{32}
   \left(a_{23}
   a_{31}-a_{21}
   a_{33}\right)}{a_{33}
   \left(a_{23}
   a_{31}+a_{21}
   a_{33}\right)}.$$
All together, there are $q\cdot (q-1)q^2\cdot q(q-1)\cdot q\cdot 1=q^5(q-1)^2$ matrices
in~$P_3'(\FF)$ which satisfy $\left(a_{23}
   a_{31}+a_{21}
   a_{33}\right)\neq 0$.

To count the matrices inside~$P_3'(\FF)$ which satisfy $\left(a_{23}
   a_{31}+a_{21}
   a_{33}\right)=0$, note that in this case $a_{21}=-\frac{a_{23}
   a_{31}}{a_{33}}$, so the first equation inside~(\ref{eq:P3'}) reduces to $\frac{2 a_{13}
   a_{23} a_{31}
   a_{32}}{a_{33}}=0$. Hence, we need to count those $3$--by--$3$ matrices which satisfy:
   $$a_{13}
   a_{23} a_{31}
   a_{32}=0,\quad a_{33}\neq 0, \quad a_{21}=-\frac{a_{23}
   a_{31}}{a_{33}},\quad a_{22}=-\frac{a_{23}
   a_{32}}{a_{33}}.$$  We may choose $(q^4-(q-1)^4)$
possible values for $a_{13},a_{23},a_{31}, a_{32}$ to have $a_{13}a_{23}a_{31} a_{32}=0$,
we may choose $(q-1)$ values for $a_{33}$, the `variables' $a_{21}, a_{22}$ are uniquely
determined,  while $a_{11}$ and $a_{12}$ are arbitrary. All together, there are
$(q^4-(q-1)^4)\cdot (q-1)\cdot 1\cdot 1\cdot q\cdot q=q^2(q-1)(q^4-(q-1)^4)$ such
matrices.

In summary we get $|P_3'(\FF)|=q^5(q-1)^2+q^2(q-1)(q^4-(q-1)^4)$. Consequently, a simple
calculation gives
\begin{align*}
|D_3(&\FF)|-|P_3(\FF)|=|D_3'(\FF)|-|P_3'(\FF)|\\
&=\Bigl(q^5(q-1)(2q-1)\Bigr)-\Bigl(q^5(q-1)^2+q^2(q-1)(q^4-(q-1)^4)\Bigr)\\
&=q^2(q-1)^5.\qedhere
\end{align*}
\end{proof}

\section{Zero permanents in $M_n(\FF)$ for $n\ge 4$ \label{Sec:3}}

The following lemma should be known, but unfortunately we were unable to find  it in the
literature. We include its proof for the sake of convenience.
\begin{lemma}\label{lem:extension-number}
Let $k\geq 2$ be an integer and let $A\in M_{k}(\FF)$ be of rank $r$. Then the set
$V^{(r)}_k(\FF)=\{({\bf x},{\bf y})\in\FF^{k}\times \FF^{k}:\;\; {\bf x}^{\tr}A{\bf
y}=0\}$ has cardinality
\begin{equation}\label{eq:Vr}
|V^{(r)}_k(\FF)|=q^{2 (k-r)}
   \left(\left(q^r-1\right)
   q^{r-1}+q^r\right).
\end{equation}
This number is a strictly decreasing function of $r$.
\end{lemma}
\begin{proof}When $r=0$, every pair satisfies ${\bf x}^{\tr}A{\bf y}=0$, so
$|V^{(r)}_k(\FF)|=q^{2k}$, which agrees with (\ref{eq:Vr}). Suppose $r>0$. There exist
invertible matrices $P$ and $Q$ such that $A=P(\Id_r \oplus 0_{k-r})Q$. Therefore,
$${\bf x}^{\tr}A{\bf y}={\bf x}^{\tr}P(\Id_r\oplus 0_{k-r})Q{\bf y}=({\bf x}')^{\tr}(\Id_r\oplus
0_{k-r}){\bf y}'$$ where ${\bf x}'=P^{\tr}{\bf x}$ and ${\bf y}'=Q{\bf y}$.
Since $P$ and $Q$ are invertible, the map $({\bf x},{\bf y})\mapsto (P^{\tr}{\bf x},Q{\bf
y})$ bijectively maps the zeros of ${\bf x}^{\tr}A{\bf y}$ onto the zeros of $({\bf
x}')^{\tr}(\Id_r\oplus 0_{n-1-r}){\bf y}'$. So we may assume that $A=(\Id_r\oplus
0_{k-r})$.

Writing ${\bf x}=(x_1,\dots,x_{k})^{\tr}$ and ${\bf y}=(y_1,\dots,y_{k})^{\tr}$ we
clearly have
$${\bf x}^{\tr}A{\bf y}=\sum_{i=1}^r x_i y_i.$$
Now, given any fixed nonzero $r$--tuple ${\bf x}_r=(x_1,\dots,x_r)$ we have that ${\bf
x}^{\tr}A{\bf y}=0$ precisely when ${\bf y}_r=(y_1,\dots,y_r)$ lies in the kernel of the
functional $F_{{\bf x}_r}:\FF^r\to\FF$ defined by ${\bf y}_r\mapsto {\bf x}_r^{\tr}{\bf
y}_r$. Hence, ${\bf y}_r$ must lie in a hyperplane inside $\FF^r$ of codimension $1$. Any
such hyperplane is isomorphic to $\FF^{r-1}$  and contains $q^{r-1}$ vectors ${\bf y}_r$.
Since there are precisely $q^r-1$ possible nonzero vectors ${\bf x}_r$, we get
$(q^r-1)q^{r-1}$ tuples $({\bf x}_r,{\bf y}_r)\in\FF^r\times \FF^r$ which satisfy
$\sum_{i=1}^r x_i y_i=0$, and such that ${\bf x}_r\ne0$. If, however, ${\bf x}_r=0$ then
${\bf y}_r$ can be arbitrary, which adds additional $1\cdot q^r$ tuples. Finally, we may
arbitrarily prescribe the values for `variables'
$x_{r+1},\dots,x_{k},\;y_{r+1},\dots,y_{k}$ giving a total of
$$((q^r-1)q^{r-1}+1\cdot q^r)\cdot q^{k-r}\cdot q^{k-r}=((q^r-1)q^{r-1}+q^r)\cdot
q^{2(k-r)}$$ tuples  $({\bf x},{\bf y})$ which solve ${\bf x}^{\tr}A{\bf y}=0$.

To prove the last statement in the lemma, we simply notice that the derivative $d/dr$ of
the above result equals $-(q-1) q^{2 k-r-1} \ln q<0$.
\end{proof}\bigskip

 We now recursively calculate the cardinality of the set $P_n(\FF) = \{A\in M_n(\FF):\, \per A=0\}$.
 Recall that the permanent can be computed with a Laplace decomposition as
\begin{equation}\label{eq:laplace}
\per A=a_{11}\per A_{11}+a_{12}\per A_{12}+\dots+a_{1n}\per A_{1n}.
\end{equation}
This suggests splitting $P_n(\FF)$ into two disjoint subsets
\begin{align*}
P_n(\FF)&=\dot{P}_n(\FF)\cup \ddot{P}_n(\FF)\\
&=\{A\in P_n(\FF):\;\;\per A_{11}\ne 0\}\cup \{A\in
P_n(\FF):\;\;\per A_{11}=0\}.
\end{align*} In $\dot{P}_n(\FF)$ we can choose
$(q^{(n-1)^2}-|P_{n-1}(\FF)|)$ blocks $A_{11}$ with nonzero permanent. For each fixed
block~$A_{11}$ we can arbitrarily prescribe the values for $2(n-1)$ `variables'
$a_{12},a_{13},\dots,a_{1n},\;a_{21},a_{31},\dots,a_{n1}$. However, the value of $a_{11}$
is then completely determined by $a_{11}=-(a_{12}\per A_{12}+\dots+a_{1n}\per
A_{1n})/\per A_{11}$. All together, the first subset has cardinality
\begin{equation}\label{eq:class1}
|\dot{P}_n(\FF)|=\big(q^{(n-1)^2}-|P_{n-1}(\FF)|\big) q^{2(n-1)}.
\end{equation}
Consider next the second set. Here, we further decompose (\ref{eq:laplace})  into
$$\per A= a_{11}\per A_{11}+\sum_{i,j=2}^n a_{1i}a_{j1}\per A_{(1i)(j1)},$$
where $A_{(1i)(j1)}$ is an $(n-2)$--by--$(n-2)$ submatrix, obtained from $A$ by deleting
the $1$-st and the $j$-th row and the $i$-th and the $1$-st column. The second term is a
bilinear form. Actually,  by introducing column vectors ${\bf
x}=(a_{12},\dots,a_{1n})^{\tr}$, ${\bf y}=(a_{21},\dots,a_{n1})^{\tr}\in\FF^{n-1}$ we can
write
$$\sum_{i,j=2}^n a_{1i}a_{j1}\per A_{(1i)(j1)}={\bf
x}^{\tr}\st{A_{11}}{\bf y},$$ where $\st{A_{11}}=\bigl(\per
A_{(1i)(j1)} \bigr)_{2\le i,j\le n}$ is the $(n-1)$--by--$(n-1)$ matrix of principal
per-minors of the matrix $A_{11}$. By Lemma~\ref{lem:extension-number} applied at
$k=n-1$, the number of pairs $({\bf x},{\bf y})$, for which the above equation is zero,
equals $|V^{(r)}_{n-1}(\FF)|=q^{2 (n-r-1)}
   \left(\left(q^r-1\right)
   q^{r-1}+q^r\right)$, where $r=\rk \st{A_{11}}$.
To count the cardinality of $\ddot{P}_n(\FF)$ we have to do the following. First, we
multiply $|V^{(r)}_{n-1}(\FF)|$ with the number of all lower-right $(n-1)$--by--$(n-1)$
blocks $A_{11}$ which have permanent equal to zero and $\rk \st{A_{11}}=r$. Then we make
a sum of these products over all ranks $r$. Finally, we multiply this sum with $q$ since
$\per A_{11}=0$ and therefore $a_{11}$ can be arbitrary. So, given an integer $r\ge 0$ we
define
\begin{equation}
{N}^{(r)}_{n-1}(\FF)=\{X\in M_{n-1}(\FF):\;\;\per X=0, \;\rk \st{X}=r\},
\end{equation}
and then
\begin{equation*}\label{eq:class2}
|\ddot{P}_n(\FF)| =q\sum_{r=0}^{n-1}   |{N}^{(r)}_{n-1}(\FF)|\cdot |V^{(r)}_{n-1}(\FF)|.
\end{equation*}
Combined with equality~(\ref{eq:class1}) for $|\dot{P}_n(\FF)|$, we  derive the following
recursive formula for number of $n$--by--$n$ matrices with permanent zero:
\begin{equation}\label{eq:Pn}
 |P_n(\FF)|=(q^{(n-1)^2}-|P_{n-1}(\FF)|) q^{2(n-1)}+q\sum_{r=0}^{n-1} |{N}^{(r)}_{n-1}(\FF)|\cdot |V^{(r)}_{n-1}(\FF)|.
\end{equation}

Unfortunately, we were unable to calculate $|P_n(\FF)|$ since we could not determine the
values for $|{N}^{(r)}_{n-1}(\FF)|$. However, in the next section we obtain an upper
bound $\mathfrak{U}_n(\FF)$ for $|P_n(\FF)|$ which is sufficient to prove the theorem.
Here is  a brief sketch of our procedure. We will introduce the functions
$\mathfrak{L}_n(\FF)$ and $\mathfrak{U}_n(\FF)$, defined inductively by the equalities
(\ref{eq:Ln}) and (\ref{eq:Un}), correspondingly, and show by the simultaneous induction
that $\mathfrak{L}_n(\FF)$ is a lower bound for $|P_n(\FF)|$ (Step 1 of the proof of
Lemma~\ref{lem:upper-lower}). Then we estimate the summand $|\dot{P}_n(\FF)|$ using the
inductively proved lower bound $\mathfrak{L}_{n-1}(\FF)$ (Step 2 of the proof of
Lemma~\ref{lem:upper-lower}). We split the summand $|\ddot{P}_n(\FF)|$ into three
summands: for $r=0$, $r=1$, and $r\ge 2$.  In order to estimate the last summand (i.e.,
with $r\ge 2$) we roughly use the monotonicity proved in Lemma~\ref{lem:extension-number}
for $|V^{(r)}_{n-1}(\FF)|$ and argue
 that a part is less than the whole, i.e., use inductive bound
$\sum_{r=2}^{n-1} |{N}^{(r)}_{n-1}(\FF)| \le |P_{n-1}(\FF)|\le \mathfrak{U}_{n-1}(\FF)$
(Step 3 of the proof of Lemma~\ref{lem:upper-lower}). Then we estimate separately the
first two summands (Steps 4 and 5 of the proof of Lemma~\ref{lem:upper-lower}).

\section{Upper and lower bounds \label{Sec:4}}
In this section we determine the upper bound $\mathfrak{U}_n(\FF)$ for $|P_n(\FF)|$. To
do this we need a lower bound $\mathfrak{L}_n(\FF)$ for $|P_n(\FF)|$ as well. To simplify
the writing we will also define auxiliary  quantities $\mathfrak{N}^{(0)}_{n-1}(\FF)$ and
$\mathfrak{N}^{(1)}_{n-1}(\FF)$. It will be shown that they are upper bounds for
$|{N}^{(0)}_{n-1}(\FF)|$ and $|{N}^{(1)}_{n-1}(\FF)|$ respectively.

First, we lower-estimate the number of $1$--by--$1$ and $2$--by--$2$ matrices with zero
permanent by $\mathfrak{L}_{1}(\FF)=0$, $\mathfrak{L}_{2}(\FF)=0$. It is easy to see that
$|P_1(\FF)|=1$ and  $|P_2(\FF)|=q^3+q^2-q$, so we define
$\mathfrak{U}_{1}(\FF)=|P_1(\FF)|=1$ and $\mathfrak{U}_{2}(\FF)=|P_2(\FF)|=q^3+q^2-q$.
Note that we have already calculated the exact value for $|P_3(\FF)|$, see
Lemma~\ref{lem:permanent_for_3-by-3} and Formula~\ref{eq:det=0}. Hence, we put
$\mathfrak{L}_{3}(\FF)=\mathfrak{U}_{3}(\FF)=|P_3(\FF)|$. We also put
$\mathfrak{N}^{(0)}_{2}(\FF)=1$. Finally, we define $\mathfrak{N}^{(0)}_{n-1}(\FF)$,
$\mathfrak{N}^{(1)}_{n-1}(\FF)$, $\mathfrak{L}_n(\FF)$, and $\mathfrak{U}_n(\FF)$ for
$n\geq 4$. We do it recursively as follows
\begin{align}
\label{eq:n_0} \mathfrak{N}^{(0)}_{n-1}(\FF)&=1+\sum_{k=1 }^{n-3}{n-1\choose
n-k-2}^2q^{2(n-k-2)(k+1)}\cdot\\
&\hspace{4cm}\mbox{}\cdot\bigl(q^{(n-k-2)^2}-\mathfrak{L}_{n-k-2}(\FF)\bigr),\notag\\[5mm]
\label{eq:n_1} \mathfrak{N}^{(1)}_{n-1}(\FF)&= \bigl(q^{(n-1)^2-1}-(q-3
   )^{(n-1)^2-1}  \bigr) +\mathfrak{N}^{(0)}_{n-2} (\FF)\cdot q^{2
   (n-2)+1}+\\[5mm]
&\hspace{5cm}\mbox{}+q\cdot\mathfrak{U}_{n-2}(\FF)\cdot |V^{(1)}_{n-2}(\FF)|,\notag\\
\label{eq:Ln} \mathfrak{L}_{n}(\FF)&=(q^{(n-1)^2}-\mathfrak{U}_{n-1}(\FF)) q^{2(n-1)},\\[5mm]
\label{eq:Un} \mathfrak{U}_n(\FF)&=(q^{(n-1)^2}-\mathfrak{L}_{n-1}(\FF)) q^{2(n-1)}+q\cdot \mathfrak{N}^{(0)}_{n-1}(\FF)\cdot |V^{(0)}_{n-1}(\FF)|+\\
\nonumber &\hspace{1cm}+q\cdot \mathfrak{N}^{(1)}_{n-1}(\FF)\cdot |V^{(1)}_{n-1}(\FF)|+q\cdot \mathfrak{U}_{n-1}(\FF)\cdot |V^{(2)}_{n-1}(\FF)|.
\end{align}
\begin{lemma}\label{lem:upper-lower}
Suppose $|\FF|>3$. Then $\mathfrak{L}_{n}(\FF)\leq |P_n(\FF)|\leq\mathfrak{U}_{n}(\FF)$
for all~$n$.
\end{lemma}
\begin{proof}
We use induction on $n$. For $n=1,2,3$ this is clear. Now,  let $n\geq 4$ and assume that
we have already proven that $\mathfrak{L}_{k}(\FF)\leq
|P_{k}(\FF)|\leq\mathfrak{U}_{k}(\FF)$ holds for all $1\le k\le n-1$. Let us show that it
holds also
 for $k=n$.

  \begin{myenumerate}{Step}
\item
 To start with, we infer from~(\ref{eq:Pn}) and from induction
hypothesis that
\begin{align*}
   |P_n(\FF)| & =|\dot{P}_n(\FF)|+|\ddot{P}_n(\FF)| \\ &  \geq |\dot{P}_n(\FF)| \\ & =  (q^{(n-1)^2}-|P_{n-1}(\FF)|) q^{2(n-1)}\\
   & \geq (q^{(n-1)^2}-\mathfrak{U}_{n-1}(\FF)) q^{2(n-1)}=\mathfrak{L}_{n}(\FF),
\end{align*}
which proves the inductive argument for the lower bound.

We now proceed with the upper bound.

 \item By the inductive hypothesis
 $$ |\dot{P}_n(\FF)| = (q^{(n-1)^2}-|P_{n-1}(\FF)|) q^{2(n-1)} \leq (q^{(n-1)^2}-\mathfrak{L}_{n-1}(\FF)) q^{2(n-1)}.$$

 \item We are using now the boundary obtained at Step 2 and split the second summand into
three parts for $r=0$, $r=1$, and $r\ge 2$ as follows
\begin{align}
\nonumber |P_n(\FF)|&\leq (q^{(n-1)^2}-\mathfrak{L}_{n-1}(\FF)) q^{2(n-1)}+q\cdot |{N}^{(0)}_{n-1}(\FF)|\cdot |V^{(0)}_{n-1}(\FF)|+\\
\nonumber &\hspace{.9cm}+q\cdot |{N}^{(1)}_{n-1}(\FF)|\cdot
|V^{(1)}_{n-1}(\FF)|+q\left(\sum_{r=2}^{n-1} |{N}^{(r)}_{n-1}(\FF)|\cdot
|V^{(r)}_{n-1}(\FF)| \right)
\end{align} Since by Lemma~\ref{lem:extension-number},
$|V^{(r)}_{n-1}(\FF)|$ is a decreasing function of $r$, we estimate its value by
$|V^{(2)}_{n-1}(\FF)|$. Since the sets ${N}^{(r)}_{n-1}(\FF)$ are obviously disjoint, we
have
$$\sum_{r=2}^{n-1} |{N}^{(r)}_{n-1}(\FF)|\le |P_{n-1}(\FF)|$$
and using the inductive bound $|P_{n-1}(\FF)|\leq \mathfrak{U}_{n-1}(\FF)$ we obtain
\begin{align}
\label{eq:up_Un} |P_n(\FF)| &\leq (q^{(n-1)^2}-\mathfrak{L}_{n-1}(\FF)) q^{2(n-1)}+q\cdot |{N}^{(0)}_{n-1}(\FF)|\cdot |V^{(0)}_{n-1}(\FF)|+\\
\nonumber &\hspace{2cm}\mbox{}+q\cdot |{N}^{(1)}_{n-1}(\FF)|\cdot |V^{(1)}_{n-1}(\FF)|+q\cdot\mathfrak{U}_{n-1}(\FF)\cdot |V^{(2)}_{n-1}(\FF)|.
\end{align}

To show  $|P_n(\FF)|\leq \mathfrak{U}_n(\FF)$ it now suffices to demonstrate that
$\mathfrak{U}_n(\FF)$, defined by~(\ref{eq:Un}), is even greater than the last quantity
in (\ref{eq:up_Un}). And to verify this  claim, it is sufficient to prove
 $|{N}^{(0)}_{n-1}(\FF)|\le \mathfrak{N}^{(0)}_{n-1}(\FF)$ and $|{N}^{(1)}_{n-1}(\FF)|\le\mathfrak{N}^{(1)}_{n-1}(\FF)$.

\item Let us prove that $|N^{(0)}_{{k}}(\FF)|\leq
\mathfrak{N}^{(0)}_{{k}}(\FF)$  for all $2\leq {k}\leq n-1$.\\

To see this, recall that $\mathfrak{L}_{k}(\FF)\leq |P_{k}(\FF)|$ for $1\le k\le n-1$ by
the inductive hypothesis. Note that $|N^{(0)}_{{k}}(\FF)|$ equals the number of all
${k}$--by--${k}$ matrices $X=\bigl(x_{ij}\bigr)$ in which every principal per-minor
vanishes. Then it is easy to see that, when ${k}=2$ all four per-minors of the
$2$--by--$2$ matrix $X$ vanish precisely when $X=0$. So, $|N^{(0)}_2(\FF)|=1$. By
definition we also have  $\mathfrak{N}^{(0)}_2(\FF)=1$. Hence, it remains to prove the
claim for $3\leq {k}\leq n-1$.

To do this we split the set $N^{(0)}_{{k}}(\FF)$ into the union of the following sets of
matrices: for any $j,\ 1\le j\le  k-2$ we consider the set of matrices with all
$({k}-i)$--by--$({k}-i)$ per-minors equal to 0 for any $i, 1\le i\le  j$ and possessing
a nonzero $({k}-j-1)$--by--$({k}-j-1)$ per-minor, and the set consisting just of the zero
matrix. Then we estimate the number of matrices in each of these sets.

We first over-estimate the number of matrices from $N_{{k}}^{(0)}(\FF)$ with the
additional property that they have a nonzero $({k}-2)$--by--$({k}-2)$ per-minor. For
simplicity assume that this $({k}-2)$--by--$({k}-2)$ submatrix is in the lower-right
corner, i.e.,
 $\per X_{(11)(22)}\ne0$; for other positions the calculations yield the same
results. Note that such $({k}-2)$--by--$({k}-2)$ lower-right block can be chosen in
$(q^{({k}-2)^2}-|P_{{k}-2}(\FF)|)$ ways. But by the inductive hypothesis, this number is
smaller or equal to $(q^{({k}-2)^2}-\mathfrak{L}_{{k}-2}(\FF))$. So, such
$({k}-2)$--by--$({k}-2)$ lower-right block can be chosen in not more than
$(q^{({k}-2)^2}-\mathfrak{L}_{{k}-2}(\FF))$ ways. By the assumption every
$({k}-1)$--by--$({k}-1)$ per-minor vanishes. In particular, $\per X_{11}=\per X_{12}=\per
X_{21}=\per X_{22}=0$, from where all the $2^2=4$ `variables' $x_{11}$, $x_{12}$,
$x_{21}$, $x_{22}$ from the upper-left $2$--by--$2$ corner are uniquely determined by the
block $X_{(11)(22)}$ and the other `variables' in the first or second row or column. For
example, $x_{22}=-\sum_{i> 2} x_{2i}\per X_{(11)(2i)} /\per X_{(11)(22)}$. Now, if we
prescribe the values for the $4({k}-2)$ `variables' $x_{i3},\dots,x_{i{k}}$ and
$x_{3i},\dots,x_{{k}i}$, $i=1,2$, arbitrarily we will obtain the estimate which is larger
or equal to the precise number. Finally, we multiply this estimate with ${{k} \choose
{k}-2}^2$ possible positions for the $({k}-2)$--by--$({k}-2)$ nonzero per-minor, to
obtain the following upper-bound:
$${{k} \choose {k}-2}^2 q^{4({k}-2)} \bigl(q^{({k}-2)^2}-\mathfrak{L}_{{k}-2}(\FF)  \bigr).$$

Among those  still remaining in our class of ${k}$--by--${k}$ matrices with all principal
per-minors zero, we next over-estimate the number of those matrices which have all
$({k}-2)$--by--$({k}-2)$ per-minors zero, but such that at least one
$({k}-3)$--by--$({k}-3)$ per-minor is nonzero. Proceeding as above, there are at most
$\bigl(q^{({k}-3)^2}-\mathfrak{L}_{{k}-3}(\FF)\bigr)$ possible such per-minors at a given
position. Having prescribed any one, there are $3^2=9$ `variables' which are completely
determined by the demand that every $({k}-2)$--by--$({k}-2)$ principal per-minor
vanishes. We may arbitrarily prescribe the values for the rest of
${k}^2-({k}-3)^2-9=6({k}-3)$ `variables.' Since there are ${{k}\choose {k}-3}^2$ possible
positions for a given nonzero $({k}-3)$--by--$({k}-3)$ per-minor, there are at most
$${{k}\choose {k}-3}^2q^{{k}^2-({k}-3)^2-9}\bigl(q^{({k}-3)^2}-\mathfrak{L}_{{k}-3}(\FF)\bigr)$$
matrices inside the present subclass. We now proceed inductively. At the ${j}$-th stage
we over-estimate those ${k}$--by--${k}$ matrices such that every $({k}-i)$--by--$({k}-i)$
per-minor vanishes, for $i=1,\dots,{j}$, but there exists  a nonzero per-minor of
dimension $({k}-{j}-1)$--by--$({k}-{j}-1)$. Arguing as above, there are at most
$${{k} \choose {k}-{j}-1}^2q^{{k}^2-({k}-{j}-1)^2-({j}+1)^2}\bigl(q^{({k}-{j}-1)^2}-\mathfrak{L}_{{k}-{j}-1}(\FF)\bigr)$$
of them. This process stops at ${j}={k}-1$, when every $1$--by--$1$ per-minor vanishes,
i.e., when $X=0$. Then we do not use the above formula because we clearly have only $1$
possibility for $X=0$. Summing up, we over-estimate $|N_{{k}}^{(0)}(\FF)|$ as
$$|N_{{k}}^{(0)}(\FF)|\leq 1+\sum_{{j}=1 }^{{k}-2}{{k} \choose
{k}-{j}-1}^2q^{2({k}-{j}-1)({j}+1)}\bigl(q^{({k}-{j}-1)^2}-\mathfrak{L}_{{k}-{j}-1}(\FF)\bigr).$$
By~(\ref{eq:n_0}) the right side equals
$\mathfrak{N}_{{k}}^{(0)}(\FF)$.\\

\item Let us prove that $|{N}^{(1)}_{n-1}(\FF)|\leq \mathfrak{N}^{(1)}_{n-1}(\FF)$.\\

To see this, we divide the set $N^{(1)}_{n-1}(\FF)$ of all $(n-1)$--by--$(n-1)$ matrices
$X$ with $\per X=0$ and $\rk\st{X}=1$ in three disjoint subsets
\begin{align*}
&\dot{N}^{(1)}_{n-1}(\FF)=\{X\in N^{(1)}_{n-1}(\FF):\;\;\per X_{11}\ne0\},\\
&\ddot{N}^{(1)}_{n-1}(\FF)=\{X\in N^{(1)}_{n-1}(\FF):\;\;\per X_{11}=0\quad \textrm{and}\quad \st{X_{11}}=0\},\\
&\dddot{N}^{(1)}_{n-1}(\FF)=\{X\in {N}^{(1)}_{n-1}(\FF):\;\;\per X_{11}=0\quad
\textrm{and}\quad \st{X_{11}}\ne 0\}
\end{align*}
and then over-estimate the cardinality of each of them.

Start with $\dot{N}^{(1)}_{n-1}(\FF)$ and recall that $\rk\st{X}\le 1$ if and only if all
$2$--by--$2$ determinant-minors of $\st{X}$ vanish. In particular, the complement of
$\dot{N}^{(1)}_{n-1}(\FF)$ inside the set ${\mathcal W}_{n-1}=\{X\in M_{n-1}(\FF):\,\per
X=0,\; \per X_{11}\ne 0\}$ contains the subset~${\mathcal V}_{n-1}$ of all
$(n-1)$--by--$(n-1)$ matrices with the following properties
\begin{align}
0&=\per X=x_{11}\per X_{11}+\sum_{i\ge 2} x_{i1}\per X_{i1}\label{eq:pom1},\\
0&\ne\per X_{11}\label{eq:pom2},\\
0&\ne\per X_{11}\per X_{22}-\per X_{12}\per X_{21}.\label{eq:pom3}
\end{align}
From~(\ref{eq:pom1})--(\ref{eq:pom2}) we express the `variable' $x_{11}$ and put it
into~(\ref{eq:pom3}). Note that the only factor in ~(\ref{eq:pom3}) which contains
$x_{11}$ is $\per X_{22}$. Therefore, after elimination of $x_{11}$ in (\ref{eq:pom3}),
the set ${\mathcal V}_{n-1}$ is determined by simultaneously non-vanishing of two
polynomials in $(n-1)^2-1$ `variables' $x_{12},\dots,x_{1(n-1)}$,
$x_{21},\dots,x_{2(n-1)},\dots\dots,x_{(n-1)(n-1)}$:
\begin{align} \label{eq:16.5}
p_1(X)=\per X_{11}&\ne0,\\
p_2(X)=\per X_{11}&(\per X_{22}|_{x_{11}=\label{eq:16.5.b}
-\frac{\sum_{i\ge 2} x_{i1}\per X_{i1}}{\per X_{11}}})\;-\\
&\hspace{3cm}\mbox{}-\;\per X_{12}\per X_{21}\ne0,\notag
\end{align}
By the definition of the permanent $p_1$ is a multilinear polynomial, i.e., every
`variable' of $p_1$ is linear, and it is also easy to see that every `variable' of $p_2$
is either linear or quadratic. Now, there exists at least one tuple of `variables' which
fulfills both inequalities. To see this, just notice that
$$X=\left(\begin{smallmatrix}
1&1\\
1&-1
\end{smallmatrix}\right)\oplus \Id_{n-3}$$
is a matrix with $\per X=0$, $\per X_{11} \ne 0$, and with $\per X_{11}\per X_{22}-\per
X_{12}\per X_{21}=-2\ne0$.

We now claim that at least $(q-3)^{(n-1)^2-1}$ tuples simultaneously satisfy both
inequalities~(\ref{eq:16.5})--(\ref{eq:16.5.b}). Namely, start with a given tuple that
does satisfy them. Keep all `variables' but one, say $x_{i_0j_0}$ for simplicity, fixed.
Recall that in the first polynomial $x_{i_0j_0}$ is at most linear, while in the second
$x_{i_0j_0}$ is at most quadratic (it may also happen that for some tuple, the
polynomials are constant). So, to satisfy the second inequality, the `variable'
$x_{i_0j_0}$ can take all but perhaps two values --- this is because a quadratic
polynomial has at most two zeros. Since the first polynomial is linear, at most one of
the allowed values of $x_{i_0j_0}$ can be its zero. So, to simultaneously satisfy also
the first inequality, we can choose at least $q-3$ values for `variable' $x_{i_0j_0}$. In
this way we obtained $(q-3)$ tuples which simultaneously satisfy inequalities
(\ref{eq:16.5})--(\ref{eq:16.5.b}).

We proceed by choosing another `variable' while keeping all the others fixed. In the same
way as before we obtain for each of the above $(q-3)$ tuples additional $(q-3)$ tuples,
hence together $(q-3)^2$ tuples which simultaneously satisfy the inequalities
(\ref{eq:16.5})--(\ref{eq:16.5.b}).

By continuing in the same manner we finally end up with at least $(q-3)^{(n-1)^2-1}$
matrices inside ${\mathcal V}_{n-1}\subseteq {\mathcal
W}_{n-1}\setminus\dot{N}^{(1)}_{n-1}(\FF)$. Recall that ${\mathcal W}_{n-1}$ is the set
of $(n-1)$--by--$(n-1)$ matrices with $\per X=0$, $\per X_{11}\ne0$. Clearly,
 $x_{11}$ is uniquely determined with the other elements of a matrix $X$, so there are at most $q^{(n-1)^2-1}$ matrices inside ${\mathcal
W}_{n-1}$. Therefore,
\begin{align*}
|\dot{N}^{(1)}_{n-1}(\FF)|&=|{\mathcal W}_{n-1}|-|{\mathcal W}_{n-1}\setminus\dot{N}^{(1)}_{n-1}(\FF)|\\
&\leq |{\mathcal W}_{n-1}|-|{\mathcal V}_{n-1}|\le q^{(n-1)^2-1}-(q-3)^{(n-1)^2-1}.
\end{align*}

 We next over-estimate the cardinality of~$\ddot{N}^{(1)}_{n-1}(\FF)$. Firstly, the number
of $(n-2)$--by--$(n-2)$ matrices $X_{11}$ with $\per X_{11}=0$ and $\st{X_{11}}=0$ equals
$|N_{n-2}^{(0)}(\FF)|$. If we enlarge such block $X_{11}$ to an $(n-1)$--by--$(n-1)$
matrix by arbitrarily prescribing the values of $2(n-2)+1$ `variables' from the first row
and column we always obtain a matrix  with permanent zero. Note that not every extension
necessarily satisfies $\rk \st{X}=1$, however we still obtain an upper bound
$|\ddot{N}^{(1)}_{n-1}(\FF)|\leq |N_{n-2}^{(0)}(\FF)|\cdot q^{2(n-2)+1}$. By Step~1,
$|N_{n-2}^{(0)}(\FF)|\leq \mathfrak{N}^{(0)}_{n-2}(\FF)$, so
$$|\ddot{N}^{(1)}_{n-1}(\FF)|\leq  \mathfrak{N}^{(0)}_{n-2}(\FF)\cdot q^{2(n-2)+1}.$$

It remains to over-estimate the cardinality of $\dddot{N}^{(1)}_{n-1}(\FF)$.  We will
make a rough  estimate. By the induction hypothesis there are at most
$\mathfrak{U}_{n-2}(\FF)$ blocks~$X_{11}$ with $\per X_{11}=0$ and $\st{X_{11}}\ne0$.
Every such block can be enlarged to $(n-1)$--by--$(n-1)$ matrix $X$ with $0=\per
X=x_{11}\per X_{11}+{\bf y}_{21}^{\tr}\st{X_{11}}{\bf x}_{12}$, by prescribing the values
for `variables' in the first row and column. Here, ${\bf y}_{21}$ is the first column
of~$X$ with the first entry removed, and ${\bf x}_{12}^{\tr}$ is the first row of $X$
with the first entry removed. Clearly then, the `variable' $x_{11}$ is arbitrary, while
the $2(n-2)$ `variables' inside ${{\bf y}_{21}}$, ${\bf x}_{12}$ must fulfill  ${\bf
y}_{21}^{\tr}\st{X_{11}}{\bf x}_{12}=0$. By the assumptions
on~$\dddot{N}^{(1)}_{n-1}(\FF)$, we have $\rk\st{X_{11}}=r\ge 1$. So, by
Lemma~\ref{lem:extension-number} there are precisely $q\cdot |V^{(r)}_{n-2}(\FF)|\le
q\cdot |V^{(1)}_{n-2}(\FF)|$ extensions. All together,
$$|\dddot{N}^{(1)}_{n-1}(\FF)|\le q\cdot\mathfrak{U}_{n-2}(\FF)\cdot |V^{(1)}_{n-2}(\FF)|,$$
wherefrom we further deduce
$$|{N}^{(1)}_{n-1}(\FF)|=|\dot{N}^{(1)}_{n-1}(\FF)|+|\ddot{N}^{(1)}_{n-1}(\FF)|+|\dddot{N}^{(1)}_{n-1}(\FF)|\leq \mathfrak{N}^{(1)}_{n-1}(\FF),$$
which ends the proof of Step~2 and consequently also the proof of the lemma. \qedhere
\end{myenumerate}
\end{proof}

\section{Proof of the main result}\label{sec-proof}

\begin{proof}[Proof of Theorem~\ref{theorem}]
By Lemma~\ref{lem:permanent_for_3-by-3}, $|P_3(\FF)|$ is strictly smaller than the number
$|D_3(\FF)|$ of $3$--by--$3$ matrices with zero determinant for arbitrary finite field
with $\charac \FF\neq 2$. This  proves the theorem in the case $n=3$. Suppose now  $n\geq
4$. Recall that $|D_n(\FF)|$ equals $q^{n^2}-\prod_{k=1}^n(q^n-q^{k-1})$. So to prove the
theorem it remains to verify that, given a fixed $n$, then for all sufficiently large $q$
one has
$$\mathfrak{U}_{n}(\FF)\lneq q^{n^2}-\prod_{k=1}^n(q^n-q^{k-1}).$$
Note that each quantity in this expression is a polynomial in $q$. It is easy to see that
\begin{equation}\label{eq:Dn}
q^{n^2}-\prod_{k=1}^n(q^n-q^{k-1})=q^{n^2-1}+q^{n^2-2}+O(q^{n^2-5}),
\end{equation}
where $O(q^k)$ is a standard notation for a quantity which satisfies\linebreak[4]
$\limsup_{q\to \infty} |O(q^k)/q^k|<\infty$. Let us prove inductively that
\begin{align*}
\mathfrak{L}_{n}(\FF)&=q^{n^2-1}-q^{n^2-2}+O(q^{n^2-3})\qquad (n\ge 4),\\
\mathfrak{U}_{n}(\FF)&=q^{n^2-1}+O(q^{n^2-3})\qquad (n\ge 4).
\end{align*}
To start with, one directly computes from~(\ref{eq:Ln}) that
\begin{align*}
\mathfrak{L}_{4}(\FF)=q^{15}-q^{14}-5 q^{12}+11
   q^{11}-9 q^{10}+4 q^9-q^8=q^{15}\!-\!q^{14}\!+\!O(q^{13})
\end{align*}
and from~(\ref{eq:Un}), (\ref{eq:n_0}), (\ref{eq:n_1}), and (\ref{eq:Vr}) that
\begin{align*}
\mathfrak{U}_{4}(\FF)&=q^{15}+53 q^{13}-520
   q^{12}+3276 q^{11}-12864
   q^{10}+\\
   &\mbox{}+32905 q^9-54445
   q^8+55410 q^7-30619
   q^6+6561 q^5\\
   &=q^{15}+O(q^{13}).
\end{align*}

Now, assume  $n\geq 5$ and the claim holds for all $\mathfrak{L}_{k}(\FF)$ and
$\mathfrak{U}_{k}(\FF)$, where $4\leq k\leq n-1$. Then,
$\mathfrak{L}_{n}(\FF)=(q^{(n-1)^2}-\mathfrak{U}_{n-1}(\FF))q^{2(n-1)}=
(q^{(n-1)^2}-q^{(n-1)^2-1}-O(q^{(n-1)^2-3}))q^{2(n-1)}=q^{n^2-1}-q^{n^2-2}+
O(q^{n^2-4})$, proving the inductive step for the lower bound.

Consider lastly $\mathfrak{U}_n(\FF)$. According to its definition~(\ref{eq:Un}), we
split it as
$$\mathfrak{U}_{n}(\FF)=I_{n}+II_{n}+III_{n}+IV_{n},$$
where
 $$I_{n}=(q^{(n-1)^2}-\mathfrak{L}_{n-1}(\FF)) q^{2(n-1)}, \quad II_{n}=q\cdot \mathfrak{N}^{(0)}_{n-1}(\FF)\cdot |V^{(0)}_{n-1}(\FF)|,$$
 $$ III_{n} (\FF)=q\cdot \mathfrak{N}^{(1)}_{n-1}(\FF)\cdot |V^{(1)}_{n-1}(\FF)|, \quad IV_n=q\cdot \mathfrak{U}_{n-1}(\FF)\cdot |V^{(2)}_{n-1}(\FF)|.$$

The first summand is done as  for $\mathfrak{L}_{n}(\FF)$ and equals
$$I_n =
q^{n^2-1}-q^{n^2-2}+ O(q^{n^2-3}).$$
 In the last summand we use~(\ref{eq:Vr}) to deduce
\begin{align*}
IV_n&=q\cdot\mathfrak{U}_{n-1}(\FF)\cdot |V^{(2)}_{n-1}(\FF)|\\
&=q\big(q^{(n-1)^2-1}+O(q^{(n-1)^2-3})\big)\cdot \big(q^{2(n-3)}(q^3+q^2-q)\big)\\
&=q^{n^2-2}+ O(q^{n^2-3}).
\end{align*}

To estimate $II_{n}$, we infer from (\ref{eq:n_0}) that
\begin{align}
\nonumber \mathfrak{N}^{(0)}_{n-1}(\FF)&=\sum_{k=1}^{n-3}O(q^{2(n-k-2)(k+1)})\cdot \big(q^{(n-k-2)^2}-O(q^{(n-k-2)^2-1})\big)\\
\label{eq:II}&=\sum_{k=1}^{n-3}O(q^{(n-1)^2-(k+1)^2})=O(q^{(n-1)^2-4}),
\end{align}
while~(\ref{eq:Vr}) implies that $|V^{(0)}_{n-1}(\FF)|=O(q^{2n-2})$. Consequently,
$II_{n}=q\cdot \mathfrak{N}^{(0)}_{n-1}(\FF)\cdot |V^{(0)}_{n-1}(\FF)|=O(q^{n^2-4})$,
which is  below the required $O(q^{n^2-3})$.

Consider lastly the third summand.  To estimate (\ref{eq:n_1}) we note that
$(q^{(n-1)^2-1}-(q-3)^{(n-1)^2-1})=O(q^{(n-1)^2-2})$. By~(\ref{eq:II}),
$\mathfrak{N}^{(0)}_{n-2}(\FF)=O(q^{(n-2)^2-4})$, while (\ref{eq:Vr}) implies
$|V^{(1)}_{n-2}(\FF)|=O(q^{2n-5})$ and $|V^{(1)}_{n-1}(\FF)|=O(q^{2n-3})$. Hence,
\begin{align*}
\mathfrak{N}^{(1)}_{n-1}(\FF)=O(q^{(n-1)^2-2})&+O(q^{(n-2)^2-4})\cdot q^{2(n-2)+1}+\\
&+q\cdot O(q^{(n-2)^2-1})\cdot O(q^{2n-5})=O(q^{(n-1)^2-2}),
\end{align*}
and
 $$III_n=q\cdot \mathfrak{N}^{(1)}_{n-1}(\FF)\cdot |V^{(1)}_{n-1}(\FF)|=q\cdot
O(q^{(n-1)^2-2})\cdot O(q^{2n-3})=O(q^{n^2-3}).$$ In total we have
$\mathfrak{U}_{n}(\FF)=I_{n}+II_{n}+III_{n}+IV_{n}=q^{n^2-1}+O(q^{n^2-3})$ which proves
the inductive step. Note that this number is strictly smaller than~(\ref{eq:Dn}) for all
sufficiently large $q$, so for such $q$ we have $|P_n(\FF)|\leq \mathfrak{U}_{n}(\FF)
\lneq |D_n(\FF)|$, which  proves the theorem.
\end{proof}

\section{Applications\label{Sec:7}}

In this section we apply the developed technique and results to estimate the probability
of the determinant and permanent functions to have a given value in a finite field.  This
problem goes back to the works of Erd\"os and R\'enyi~\cite{erdos_renyi,erdos_renyi1},
where they estimated the probability for a (0,1)-matrix with a given number of ones to
have a zero permanent. Later many authors investigated this topic for determinant and
permanent functions of (0,1)-matrices, see monographs~\cite{borovskikh_korolyuk,girko}
for details. In particular, Sachkov~\cite{sachkov} proved that if a uniform distribution
is given on the set of all (0,1)-matrices of size $m\times n$, where $m\le n$, then the
probability $P\{\per A\ne 0\} \to 1 $ if $n\to \infty$, where $A$ is an arbitrary
(0,1)-matrix of size $m\times n$. An asymptotics for cardinality of (0,1)-matrices with
zero permanent was given by Everett and  Stein in~\cite{everett_stein}, corresponding
results for the determinant are due to Koml\'os, see~\cite{komlos,komlos2}.

 Here  we investigate the situation over arbitrary finite fields. The
application of our technique over a finite field $\FF$ of cardinality $q$ shows that for
$0\ne \alpha\in \FF$ the probability function $P$ behaves as follows
$$P(\det A=\alpha)= \frac{1}{q}-\frac{1}{q^3}+O(\frac{1}{q^4}), \qquad P(\det A=0)= \frac{1}{q}+\frac{1}{q^2}+O(\frac{1}{q^5}),$$
$$  \frac{1}{q}-\frac{1}{q^2}+O(\frac{1}{q^3}) \le  P(\per A=0) \le \frac{1}{q}+O(\frac{1}{q^3}), \mbox{ and } $$
$$ \frac{1}{q}+O(\frac{1}{q^4}) \le P(\per A=\alpha)\le \frac{1}{q}+\frac{1}{q^3}+O(\frac{1}{q^4}),$$
so, roughly speaking, each of these probabilities approximately equals to~$1/q
+O(\frac{1}{q^2})$.

In order to prove our result we need the following lemma:
\begin{lemma} \label{all_eq}
Let $\FF$ be a finite field, $\charac \FF\ne 2$. Then for any nonzero $\alpha,\beta\in
\FF$ the cardinality of the set of matrices of a given size with the determinant
(permanent) $\alpha$ is equal to the cardinality of the set of matrices of a given size
with the determinant (permanent) $\beta$, i.e.,
$$
|\{A\in M_n(\FF):\;\; \det A=\alpha\}|=|\{A\in M_n(\FF):\;\; \det A=\beta\}|$$
and
$$
|\{A\in M_n(\FF):\;\; \per A=\alpha\}|=|\{A\in M_n(\FF):\;\; \per A=\beta\}|.$$
\end{lemma}

\begin{proof}
We denote $D^\alpha_n(\FF)=\{A \in M_n(\FF):\;\; \det A=\alpha\}$.

For any $A=(a_{ij})\in D^\alpha_n(\FF)$ we consider the matrix $B=(b_{ij})$ defined by
$b_{ij}=a_{ij}$ for $i=1,\ldots, n$, $j=2,\ldots, n$, $b_{i1}=\frac{\beta}{\alpha}a_{i1}$
for $i=1,\ldots, n$. Then $\det B=\frac{\beta}{\alpha} \det A=\beta$, i.e., $B\in
D^\beta_n(\FF)$. Since $\alpha\beta\ne 0$, the mapping from $A$ to $B$ is well-defined
and injective, hence, $|D^\beta_n(\FF)|\ge |D^\alpha_n(\FF)|$. Similarly,
$|D^\beta_n(\FF)|\le |D^\alpha_n(\FF)|$.

Since permanent is also a linear function of a row or a column, the result for the
permanent can be obtained in the same way.
\end{proof}

\begin{theorem}
Let $\FF$ be a finite field, $|\FF|=q$, $\charac \FF\ne 2$. For any $\alpha\in \FF$ the
probability that $\det A=\alpha$, $A\in M_n(\FF)$, is equal to $\frac{1}{q} +
O(\frac{1}{q^2})$ and the probability that $\per A=\alpha$ is also equal to $\frac{1}{q}
+ O(\frac{1}{q^2})$.
\end{theorem}

\begin{proof}
We consider at first $\alpha=0$. Then by the proof of Theorem~\ref{theorem} it follows
that the quantity of matrices with zero determinant
$|D_n(\FF)|=q^{n^2-1}+q^{n^2-2}+O(q^{n^2-5})$. Hence, the probability
\begin{align*}
P(\det A=0)&= \frac{q^{n^2-1}+q^{n^2-2}+O(q^{n^2-5})}{q^{n^2}}\\
&=\frac{1}{q}+\frac{1}{q^2}+O(\frac{1}{q^5})\\
&=\frac{1}{q}+O(\frac{1}{q^2}).
\end{align*}

Similarly, using the proof of Theorem~\ref{theorem} we have
\begin{align*}
P(\per A=0)&\le \frac{\mathfrak{U}_{n}(\FF)}{q^{n^2}}= \frac{q^{n^2-1}+O(q^{n^2-3})}{q^{n^2}}=\\
&\le\frac{1}{q}+O(\frac{1}{q^3})=\frac{1}{q}+O(\frac{1}{q^2})
\end{align*}
and
\begin{align*}
P(\per A=0)&\ge \frac{\mathfrak{L}_{n}(\FF)}{q^{n^2}}= \frac{q^{n^2-1}
-q^{n^2-2}+O(q^{n^2-3})}{q^{n^2}}=
\frac{1}{q}-\frac1{q^2}+O(\frac{1}{q^3})\\
&\ge\frac{1}{q}+O(\frac{1}{q^2}).
\end{align*}
So, $P(\per A=0)= \displaystyle \frac{1}{q}+O(\frac{1}{q^2}).$

If $\alpha\ne 0$ then by Lemma~\ref{all_eq}
\begin{align*}
|D_n^\alpha(\FF)|&=\displaystyle
\frac{\prod_{k=1}^n(q^n-q^{k-1})}{q-1}= q^{\frac{n(n-1)}{2}} (q^n-1)\cdots (q^2-1)\\
&=q^{n^2-1}-q^{n^2-3}+O(q^{n^2-4}).
\end{align*} Thus  the probability
\begin{align*}
P(\det A=\alpha)&= \frac{q^{n^2-1}-q^{n^2-3}+O(q^{n^2-4})}{q^{n^2}}=\frac{1}{q}-\frac{1}{q^3}+O(\frac{1}{q^4})\\
&=\frac{1}{q}+O(\frac{1}{q^2}).
\end{align*}
Finally,
$$|P_n^\alpha(\FF)|\le \displaystyle \frac{q^{n^2}-q^{n^2-1}+q^{n^2-2}+O(q^{n^2-3}
)}{q-1}= q^{n^2-1}+q^{n^2-3}+O(q^{n^2-4})$$
and
$$|P_n^\alpha(\FF)|\ge \displaystyle \frac{q^{n^2}-q^{n^2-1}+O(q^{n^2-3})}{q-1}= q^{n^2-1}+O(q^{n^2-
4}).$$

Thus  the probability
$$P(\per A=\alpha)= \frac{1}{q}+O(\frac{1}{q^2}).$$
\end{proof}

\section{Examples and Remarks\label{Sec:6}}

\begin{remark}
In the table below, for a given $n$ we  compute the first integer $i$ such that for any
$j>i$ the value of the polynomial $ \mathfrak{U}_{n}(\FF)$ at $q=j$ is strictly less than
the value of the polynomial $|D_n(\FF)|$ at $q=j$. In the third row we give the minimal
number of elements in a field with this property, i.e., the minimal power of a prime
$q=|\FF|$ such that $\mathfrak{U}_n(\FF)\lneq |D_n(\FF)|$. We used Wolfram's
\emph{Mathematica 5.1} for the calculations. For example, when $n=5$ we have
$\mathfrak{U}_5(\FF)\lneq |D_5(\FF)|$  whenever the finite field $\FF$ has at least $76$
elements and its characteristic differs from $2$. The smallest such field with at least
$76$ elements is $GF(79)$. So, $q=79$. {\tiny$$\begin{array}{r|cccccccccccccccccccc}
 n & 3& 4 & 5 &6  &7  &8  &9  &10  &11 \\\hline
 i & 2& 43& 76&116&164&221&287&362 &446        \\\hline
 q & 3& 43& 79&121&167&223&289&367 &449        \\\hline\hline
n&12 &13&14&15&16&17&18&19&20\\\hline i&538&640&750&869&996&1133&1278&1433&1596
\\\hline
 q & 541& 641 & 751& 877&997&1151&1279&1433&1597
\end{array}$$}
\end{remark}

If $\FF$ is an infinite field then there do exist bijective converters of permanent into
determinant. In the Example~\ref{Ex:main} we give such bijective maps $\Phi: M_n(\FF)\to
M_n(\FF)$, $n\ge 2$, that even satisfy $\per A=\det \Phi(A)$ and $\det A=\per \Phi(A)$
simultaneously for all $A\in M_n(\FF)$.

\begin{example} \label{Ex:main}
If $\charac\FF = 2$ then $\per A=\det A$ for any $A\in M_n(\FF)$ so we take $\Phi(X)=X$
to achieve $\per A=\det \Phi(A)$ and $\det A=\per \Phi(A)$. Assume $\charac\FF\ne2$. Note
that the cardinality of infinite sets satisfies $|\FF|=|\FF\times\FF|$, so
$|M_n(\FF)|=|\FF^{n^2}|=|\FF|$.

We are going to prove now that for any given $\lambda,\mu\in \FF$ the cardinality of the
set of matrices with permanent $\lambda$ and determinant $\mu$ is equal to $|\FF|$, so
for any given pair of such sets there is a bijection between them. Let us denote
$$\Omega_n(\lambda,\mu)=\{A\in M_n(\FF):\;\; \per(A)=\lambda \mbox{ and } \det
(A)=\mu\}.$$

1. For given fixed $\lambda,\mu\in\FF$ consider the set
$$\Delta_n(\lambda,\mu)=\left\{\left(\begin{smallmatrix}
\alpha &(\lambda-\mu)/2\\
1 &(\lambda+\mu)/(2\alpha)
\end{smallmatrix}\right)\oplus \Id_{n-2}:\;\;\alpha\in\FF\backslash\{0\}\right\}\subseteq M_n(\FF).$$

 2. The cardinality of this set is $|\FF|-1=|\FF|$ and every matrix
from this set has permanent and determinant equal to $\lambda$ and $\mu$ respectively.

 3. Moreover,
$$\Delta_n(\lambda,\mu)\subseteq\Omega_n(\lambda,\mu)\subseteq M_n(\FF),$$ and comparing cardinalities, we obtain
$|\FF|=|\Delta_n(\lambda,\mu)|\preceq|\Omega_n(\lambda,\mu)| \preceq |M_n(\FF)|=|\FF|$.
By the classic Bernstein-Schroeder's theorem~\cite[Cor. II.7.7]{dugundji} we have
$|\Omega_n(\lambda,\mu)|=|\FF|$.

 4. So, there is a bijection
$\Phi_{\lambda,\mu}: \Omega_n(\lambda,\mu) \to \Omega_n(\mu,\lambda)$.

 5. However, due to partition
$$M_n(\FF)=\bigcup_{\lambda,\mu\in\FF} \Omega_n(\lambda,\mu),$$
the maps $\Phi_{\lambda,\mu}$ constitute a well-defined bijection $\Phi:M_n(\FF)\to
M_n(\FF)$ with $\per A=\det\Phi(A)$ and $\det A=\per\Phi(A)$. It is given by
$A\mapsto\Phi_{\lambda,\mu}(A)$ if $A$ satisfies $\per(A)=\lambda$ and $\det(A)=\mu$.
\end{example}

\begin{remark}
By adopting the above arguments it can be shown that there exists a bijection
$\Phi:M_n(\FF)\to M_m(\FF')$ with similar properties as in the previous example, provided
that $\FF$ and $\FF'$ are infinite fields of the same cardinality and $m,n\ge 2$.
\end{remark}
 Note that for any field $\FF$ there exist nonbijective converters
of permanent into determinant.
\begin{example} \label{Ex:1}
As an example, $\Phi: A\mapsto(\Id_{n-1}\oplus\per A)$ satisfies $\per A=\det\Phi(A)$.
Note that such transformations cannot be linear.
\end{example}

 Moreover,  there exist also nonbijective transformations $\Phi:
M_n(\FF)\to M_m(\FF)$  which exchange permanent and determinant of a matrix.
\begin{example} \label{Ex:2}
 If $\charac\FF = 2$ then $\per A=\det A$ for any $A\in M_n(\FF)$
and  if $m\ge n$  the map $\Phi: A \mapsto A\oplus I_{m-n}$ has the required property.

If $\charac\FF\ne2$ then for any field $\FF$ and for all $m\ge 2$ we consider
$$\Phi: A\mapsto \left(\begin{array}{cc} 1 & \frac{1}{2}(\det A - \per A) \\[1mm] 1 & \frac{1}{2}(\det A + \per A) \end{array} \right) \oplus \Id_{m-2}.$$ Hence, $\Phi$ satisfies $\per A=\det\Phi(A)$ and $\det A=\per\Phi(A)$. Note that such
transformations cannot be linear.
\end{example}

In order to extend our results to finite rings we need the following lemma, which we
include here with its proof for the sake of completeness.
\begin{lemma} \label{lem_rings}
Let $R$ be a finite ring without zero divisors. Then $R$ is a field.
\end{lemma}
\begin{proof}
Since $R$ has no zero divisors, then for any $a\in R$, $a\ne 0$, the transformations
$r_a: x\to ax$ and $l_a: x\to xa$ are injective. Thus both these transformations are
bijective since they are surjective by the finitness of $R$.

Let us check that  the neutral element  is automatically in $R$. Since $r_a$ is
surjective, there exists $x\in R$ such that $ax=a$. Now, for any $b\in R$ there exists
$y\in R$ such that $b=ya$. Thus $bx=yax=ya=b$, i.e., $x$ is a right unity. Similarly,
there is $x'\in R$ which is a left unity. Then $x=x'x=x'$, i.e., $x$ is a unity. Let us
denote it by $e$.

Now for any $a\in R$, $a\ne 0$, there exist $a', a''\in R$, such that $aa'=e$ and
$a''a=e$ by the surjectivity of $r_a$ and $l_a$, correspondingly. Considering
$a''=a''(aa')=(a''a)a'=a'$, we get that $a$ is invertible. Thus $R$ is a division ring.
By Wedderburn theorem any finite division ring is a field and the result follows.
\end{proof}

\begin{proof}[Proof of Corollary~\ref{cor_rings}]
It follows directly by the application of Theorem~\ref{theorem} to the result of
Lemma~\ref{lem_rings}.
\end{proof}

\begin{remark}
By Lemma~\ref{lem_rings} the results of Section~\ref{Sec:7} are valid for finite rings
without zero divisors as well.
\end{remark}

\end{document}